\documentclass[10pt]{amsart}
\usepackage{mathptmx}
\usepackage{amsmath}     
\usepackage{amssymb}

\begin{document}

\title{Upper bounds and asymptotic expansion for  Macdonald's  function and the summability of the Kontorovich-Lebedev integrals}

\author{ {\bf S. Yakubovich} \\
 {\em {Department of Mathematics, Fac. Sciences of University of Porto, Rua do Campo Alegre,  687; 4169-007 Porto (Portugal) }}}
 \thanks{E-mail: syakubov@fc.up.pt} 
\vspace{3mm}

\begin{abstract} {\noindent Uniform upper bounds and the asymptotic expansion with an explicit remainder term are established for the Macdonald function $K_{i\tau}(x)$. The results can be  applied, for instance,  to study the summability of the divergent Kontorovich-Lebedev integrals in the sense of Jones. Namely,  we answer affirmatively a  question (cf. [6])   whether these integrals converge for even entire functions of the exponential type in a weak sense.}
\end{abstract}

\vspace{4mm}

\maketitle

\markboth{\rm \centerline{ S. Yakubovich}}{ Kontorovich-Lebedev integrals}

{\bf Keywords}:  { Kontorovich-Lebedev transform, modified Bessel function, Macdonald function, Jones summability}

{\bf AMS Subject Classification}:   {44A15,  41A60,  33C10, 33C15}

\vspace{4mm}

\section{Upper bounds of the Lebedev type}

In the theory of the Kontorovich-Lebedev transform (KL-transform) [1], whose inverse involves the integration with respect to the pure imaginary index of the modified Bessel function or Macdonald function $K_{i\tau}(x)$ it is  important to  know its  asymptotic behavior at infinity by the index for a bounded positive argument,  or  to have an upper bound, being valued for  arbitrary real parameter $\tau$  and positive argument $x$. For instance, to investigate mapping properties  in Lebesgue spaces of the Kontorovich-Lebedev operator   we oftenly use the so-called Lebedev inequality (cf. [1], p. 219)

$$\left|K_{i\tau}(x)\right| \le A {x^{-1/4} \over \sqrt {\sinh(\pi \tau)}},\quad x, \tau > 0,\eqno(1.1)$$
where $A >0$ is the absolute constant which will be defined below.  As far as the author is aware,  such bounds are unknown for the Mehler-Fock, index Whittaker, Olevskii transforms, Lebedev's transform with the product of the modified Bessel functions among the others index  transformations [1], and will be established in the forthcoming paper.  Moreover, Lebedev conjectured in 1951 in his Habilitation Thesis more tight  than  (1.1) bounds by the index 

$$\left|K_{i\tau}(x)\right| \le A {x^{\mp 1/4} \over \sqrt {\tau \sinh(\pi \tau)}},\quad 0 < x \le 1\  (x  \ge 1), \ \tau > 0.\eqno(1.2)$$
Our goal will be to establish such kind of  upper bounds being valued for arbitrary positive parameters in the index and argument.   To do this,  we will appeal to the bounds for Bessel function of the first kind [2], namely  to estimate the value

$$c_\nu= \sup_{x \ge 0} \sqrt x \left|J_\nu(x)\right|, \quad \nu \ge - {1\over 2}.\eqno(1.2)$$ 
The classical Szeg{\" o} result says that $c_\nu= \sqrt{2/\pi}$ for $|\nu| \le 1/2$.   For $\nu >0$ Olenko proved the following inequality

$$c_\nu \le  b \left( \nu^{1/3} + {\alpha\over \nu^{1/3} } + {3 \alpha^2\over 10 \nu }\right)^{1/2},\eqno(1.3)$$
where $b= 0.674885\dots $ is the  Landau constant and $\alpha= 1.855757\dots$\ .    

Let us give an explicit value for $A$ in the Lebedev inequality (1.1). To proceed this, we will  consider the following integral (cf. [3], Vol. II, Entry 2.12.14.2)

$$K^2_{i\tau}(x) = {\pi\over \sinh(\pi\tau)} \int_0^\infty J_0(2x \sinh t) \sin(2 \tau t) dt.\eqno(1.4)$$
Hence via elementary substitutions in the integral (1.4), and since the modified Bessel function $K_{i\tau}(x),\ x, \tau > 0$ is real-valued,   we have the estimate

$$K^2_{i\tau}(x) \le  {\pi\over \sinh(\pi\tau)} \int_0^\infty \left|J_0(2x \sinh t)\right|  dt =  {\pi\over \sinh(\pi\tau)} \int_0^\infty \left|J_0(2x y)\right|  {dy\over (y^2+1)^{1/2} }$$

$$\le  {\pi \ \sup_{t\ge 0} \sqrt t \left|J_0(t)\right| \over (2x)^{1/2} \sinh(\pi\tau)} \int_0^\infty  {dy\over y^{1/2} (y^2+1)^{1/2} }.$$
The latter supremum is  $c_0= \sqrt{2/\pi}$ via (1.2). Therefore, calculating a simple beta-integral (see [3], Vol. I, Entry 2.2.4.24), we derive the  Lebedev inequality (1.1) with an explicit constant 

 $$\left|K_{i\tau}(x)\right| \le  {\Gamma(1/4)\over \sqrt{2}}\  {x^{-1/4} \over \sqrt{\sinh (\pi\tau)}},\quad x, \tau > 0,\eqno(1.5)$$
where $\Gamma(z)$ is the Euler gamma function [3], Vol. III.  Now we will derive a family of Lebedev's upper bounds,  appealing to the integral in [3], Vol. II, Entry 2.12.4.28, namely

$$K_{\nu-\rho+1}(x)= \left({2\over x}\right)^{\rho-1} \Gamma(\rho)  \int_0^\infty {y^{\nu+1}\  J_\nu(xy) \over (1+y^2)^{\rho}} dy,\eqno(1.6)$$
where $x >0, \ -1 <  \nu < 2  {\rm Re} \rho - 1/2$.  In fact, letting $\rho = \nu - \mu +1 + i\tau,\ \mu < (\nu +1/2)/2$, we have from (1.3), (1.6)

$$\left| K_{\mu+i\tau}(x) \right| \le 2^{\nu-\mu}  c_\nu  \bigg| \Gamma\left(\nu - \mu +1 + i\tau \right)\bigg|  x^{\mu-\nu-1/2} \int_0^\infty {y^{\nu+1/2} \over (1+y^2)^{\nu-\mu+1}} dy $$

$$=  2^{\nu-\mu-1}  c_\nu {\Gamma((\nu+3/2)/2) \Gamma((\nu- 2\mu+ 1/2)/2)\over \Gamma(\nu-\mu+1)} \bigg| \Gamma\left(\nu - \mu +1 + i\tau \right)\bigg|  x^{\mu-\nu-1/2},$$
i.e. we derive the following inequality

$$ \left| K_{\mu+i\tau}(x) \right| \le 2^{\nu-\mu-1}  c_\nu {\Gamma((\nu+3/2)/2) \Gamma((\nu-2\mu+ 1/2)/2)\over \Gamma(\nu-\mu+1)} $$

$$\times \bigg| \Gamma\bigg(\nu - \mu +1 + i\tau \bigg)\bigg|  x^{\mu-\nu-1/2},\quad x, \tau > 0.\eqno(1.7)$$
When $\tau \to +\infty$ the gamma function in (1.7) behaves as $O( \tau^{\nu-\mu+1/2} \exp (-\pi\tau/2))$. Consequently, letting $\mu=0$ in (1.7),  and taking into account values of $c_\nu$ in (1.3), it gives after a slight simplification  the following upper bounds

$$  \left| K_{i\tau}(x) \right| \le { \Gamma(\nu+ 1/2)\over \Gamma(\nu+1)}\  \bigg| \Gamma\left(\nu  +1 + i\tau \right)\bigg|  x^{-\nu-1/2},\   - {1\over 2}  < \nu \le {1\over 2},\eqno(1.8)$$

$$  \left| K_{i\tau}(x) \right| \le   b \left( \nu^{1/3} + {\alpha\over \nu^{1/3} } + {3 \alpha^2\over 10 \nu }\right)^{1/2} { \Gamma(\nu+ 1/2)\over \Gamma(\nu+1)}$$

$$\times   \bigg| \Gamma\left(\nu  +1 + i\tau \right)\bigg|  x^{-\nu-1/2},\quad   \nu > 0.\eqno(1.9)$$ 
If $\nu =0$ in (1.8), it yields the modified Lebedev inequality

$$   \left| K_{i\tau}(x) \right| \le \sqrt {{\pi^2 \tau\over  x \sinh(\pi\tau)}},\quad x, \tau > 0.\eqno(1.10)$$
Let $\nu= -1/2 + \delta,\  0 < \delta  < 1$.  Hence, inequality (1.8) implies

$$  \left| K_{i\tau}(x) \right| \le {  \Gamma(\delta) \  x^{- \delta}  \over \Gamma(\delta+ 1/2)} \bigg| \Gamma\left({1\over 2} +\delta + i\tau \right)\bigg| , \quad x, \tau > 0.\eqno(1.11)$$
Let $\mu=\nu < 1/2$.  Then we find from (1.7)

$$  \left| K_{\nu+i\tau}(x) \right| \le \Gamma((\nu+3/2)/2) \Gamma((1/2-\nu)/2)\  \sqrt {{ \tau\over  2 x \sinh(\pi\tau)}},\quad x, \tau > 0.\eqno(1.12)$$
One can put $\nu= -1/2$ in (1.12). This yields immediately uniform upper bounds for the kernels ${\rm Re} K_{1/2+i\tau}(x),\  {\rm Im} K_{1/2+i\tau}(x)$ of the Lebedev-Skalskaya transforms [1]. Namely, we have, for instance for the ${\rm Re}$-transform

$$\left|{\rm Re} K_{1/2+i\tau}(x)\right| = \left|{\rm Re} K_{-1/2-i\tau}(x)\right| = \left|{\rm Re} K_{-1/2+i\tau}(x)\right|$$

$$ \le \left| K_{-1/2+i\tau}(x) \right| \le   \sqrt {{ \pi^2 \tau\over  2 x \sinh(\pi\tau)}},$$
i.e. (as well as  for the ${\rm Im}$-transform )

$$ \left|{\rm Re} K_{1/2+i\tau}(x)\right|  \le   \sqrt {{ \pi^2 \tau\over  2 x \sinh(\pi\tau)}}, \quad x, \tau > 0,\eqno(1.13)$$

$$ \left|{\rm Im} K_{1/2+i\tau}(x)\right|  \le   \sqrt {{ \pi^2 \tau\over  2 x \sinh(\pi\tau)}}, \quad x, \tau > 0.\eqno(1.14)$$
The case $\mu=1,\ \nu > 3/2 $ in (1.7) gives the inequality

$$ \left| K_{1+i\tau}(x) \right| \le 2^{\nu- 2}  c_\nu {\Gamma((\nu+3/2)/2) \Gamma((\nu- 3/2)/2)\over \Gamma(\nu)}  \bigg| \Gamma\bigg(\nu  + i\tau \bigg)\bigg|  x^{1/2-\nu}.\eqno(1.15)$$
Taking into account properties for the modified Bessel functions [2], Vol. II, namely, the equality

$$\tau  K_{i\tau} \left(x\right) = x \ {\rm Im}\ K_{1+i\tau} \left(x\right),\eqno(1.16)$$
we obtain from (1.15) 

$$  \left| K_{i\tau}(x) \right| \le 2^{\nu- 2}  c_\nu {\Gamma((\nu+3/2)/2) \Gamma((\nu- 3/2)/2)\over \tau \Gamma(\nu)}  \bigg| \Gamma\bigg(\nu  + i\tau \bigg)\bigg|  x^{3/2-\nu}.\eqno(1.17)$$
If $\nu= 3/2 +\delta,\ \delta > 0$, we have (compare with (1.11))

$$   \left| K_{i\tau}(x) \right| \le 2^{-1/2+\delta}  c_\nu {\Gamma((3+\delta)/2 ) \Gamma(\delta/2)\over \tau \Gamma(3/2+\delta)}  \bigg| \Gamma\bigg(3/2+\delta  + i\tau \bigg)\bigg|  x^{-\delta}, \  x, \tau > 0.\eqno(1.18)$$

Further,  returning to the unproved Lebedev bound (1.2),  we recall integral (1.4) and the derivative of the Bessel function $J_0(z)$ to write it in terms of the iterated improper integral

$$ {1\over \pi} \sinh(\pi\tau) K^2_{i\tau}(x) = -  \lim_{T\to \infty} \int_0^T  \sin(2 \tau t) \int_{0}^{2x \sinh t} J_1(y)  dy dt.\eqno(1.19)$$
Interchanging the order of integration in (1.19) owing to the uniform convergence by $t \in [0,T]$, it reads

$$ {1\over \pi} \sinh(\pi\tau) K^2_{i\tau}(x) =  \lim_{T\to \infty} \int_0^T  \sin(2 \tau t) \left[ 1- \int_{0}^{2x \sinh t} J_1(y)  dy \right] dt$$

$$=   \lim_{T\to \infty} \bigg[ {1- \cos(2\tau T) \over 2\tau} - \int_{0}^{T} J_1(y) \int_{\log(y/(2x)+ (y^2/ (4x^2) +1)^{1/2})}^T  \sin(2 \tau t) dt dy \bigg].$$
Hence the integration with respect to $t$  and simple change of variables yield

$$ {1\over \pi} \sinh(\pi\tau) K^2_{i\tau}(x) =    \lim_{T\to \infty} \bigg[ {1- \cos(2\tau T) \over 2\tau} +  {x\over \tau} \int_{0}^{T/(2x)} J_1(2xy) \bigg[ \cos \left(2 \tau T\right) \bigg.\bigg.$$

$$\bigg. \bigg. - \cos \left(2 \tau \log \left( y + \left( y^2  +1\right)^{1/2}\right)\right) \bigg] dy\bigg] =   \lim_{T\to \infty}  {1\over 2\tau} \bigg[ 1 -   \cos \left(2 \tau T\right)   J_0\left(T\right) $$

$$\bigg. - 2x \int_{0}^{T/(2x)}  J_1(2xy)  \cos \left(2 \tau \log \left( y + \left( y^2  +1\right)^{1/2}\right)\right)  dy\bigg].\eqno(1.20)$$
Passing to the limit, we find the following representation

$$  {\tau\over \pi} \sinh(\pi\tau) K^2_{i\tau}(x) = 1 -  2x \int_{0}^{\infty}  J_1(2xy) $$

$$\times \cos \left(2 \tau \log \left( y + \left( y^2  +1\right)^{1/2}\right)\right)  dy.\eqno(1.21)$$
In order to proceed further estimations we will employ the equality for Bessel function from [4]

$$J_\nu(z)= \left({2\over \pi z}\right)^{1/2} \left( \cos \left( z - {\pi\over 2} \nu - {\pi\over 4} \right) \left( \sum_{n=0}^{N-1} (-1)^n {a_{2n}(\nu)\over z^{2n}} + R_{2N} (z,\nu)\right)\right.$$

$$\left. -  \sin \left( z - {\pi\over 2} \nu - {\pi\over 4} \right) \left( \sum_{n=0}^{M-1} (-1)^n {a_{2n+1}(\nu)\over z^{2n+1}} - R_{2M+1} (z,\nu)\right)\right),\eqno(1.22)$$
where 

$$a_n(\nu)= (-1)^n {\cos(\pi\nu)\over 2^n n! \pi} \Gamma\left(n+{1\over 2}+\nu\right) \Gamma\left(n+{1\over 2}-\nu\right),\eqno(1.23)$$
and the remainder is given by the integral [4, formula (23)]

$$R_N(z,\nu)= (-1)^{[N/2]} \left({2\over \pi }\right)^{1/2}  {\cos(\pi\nu)\over z^N \pi} \int_0^\infty {t^{N-1/2} e^{-t} \over 1+ (t/z)^2} K_\nu(t) dt,\eqno(1.24)$$
provided $|{\rm Re} \nu | < N- 1/2,\  |\arg z| < \pi/ 2$.  Therefore we write by virtue of (1.3), (1.21), (1.22), (1.24)

$$  {\tau\over \pi} \sinh(\pi\tau) K^2_{i\tau}(x) \le 1 +  2 b \sqrt{2x}  \left( 1 + \alpha + {3 \alpha^2\over 10}\right)^{1/2}$$

$$+   \left({4 x\over \pi }\right)^{1/2}  \left| \int_{1}^{\infty}    \sin \left( 2xy -  {\pi\over 4} \right) \cos \left(2 \tau \log \left( y + \left( y^2  +1\right)^{1/2}\right)\right)  {dy \over \sqrt y} \right| $$

$$+  \left({4 x\over \pi }\right)^{1/2}  \left[ \sum_{n=1}^{N-1}  {\Gamma\left(2n+{3\over 2}\right) \Gamma\left(2n-{1\over 2}\right) \over 4^n (2n)! (2n-1/2)} + \sum_{n=0}^{M-1}  {\Gamma\left(2n+{5\over 2}\right) \Gamma\left(2n+{1\over 2}\right) \over 4^{n}  (2n+1)! (4n+1)}\right]$$

$$+   {2\sqrt {2x} \over \pi^2} \int_{1}^{\infty}  \int_0^\infty {y\  e^{-t} K_1(t)  \over   y^2+ t^2}\ \left[ \left({t\over y}\right)^{2N-1/2} + \left({t\over y}\right)^{2M+1/2} \right] dt dy,$$
where $x,\ \tau > 0,\  M, N \in \mathbb{N}$.  The latter double integral can be estimated via Entry 8.4.23.3 in [3], Vol. III and an elementary inequality. Thus we get finally

$$   {\tau\over \pi} \sinh(\pi\tau) K^2_{i\tau}(x) \le 1 +  2 b \sqrt{2x}  \left( 1 + \alpha + {3 \alpha^2\over 10}\right)^{1/2}$$

$$+  {4\sqrt x \over \pi } \left[ \sum_{n=1}^{N-1}  2^{2n-1} \Gamma \left(2n-{1\over 2}\right) B\left(2n+{3\over 2}, 2n- {1\over 2}\right) \right.$$

$$\left.  +  \sum_{n=0}^{M-1}   4^n  \Gamma\left(2n+{1\over 2}\right) \ B\left(2n+{5\over 2}, 2n+{1\over 2}\right)\right]$$

$$+   {2\sqrt {x} \over \pi^2 }\left[ 4^N \  \Gamma \left(2N-{3\over 2}\right) B\left(2N-{3\over 2}, \ 2N+{1\over 2}\right) \right. $$

$$\left. +   4^{M-1} \Gamma \left(2M-{1\over 2}\right) B\left(2M-{1\over 2}, \ 2M+{3\over 2}\right)\right]$$

$$+   \left({4 x\over \pi }\right)^{1/2}  \left| \int_{1}^{\infty}    \sin \left( 2xy -  {\pi\over 4} \right) \cos \left(2 \tau \log \left( y + \left( y^2  +1\right)^{1/2}\right)\right)  {dy \over \sqrt y} \right|,\eqno(1.25) $$
where $B(z,w)$ is the Euler beta function [3]. The main goal is to estimate the integral in (1.25). To do this, we  fix a positive $\delta$ an split the integral as follows

$$ \int_{1}^{\infty}    \sin \left( 2xy -  {\pi\over 4} \right) \cos \left(2 \tau \log \left( y + \left( y^2  +1\right)^{1/2}\right)\right)  {dy \over \sqrt y} $$

$$= \left( \int_{1}^{(1+\tau)^{2\delta}} +  \int_{(1+\tau)^{2\delta}}^\infty  \right)  \sin \left( 2xy -  {\pi\over 4} \right) \cos \left(2 \tau \log \left( y + \left( y^2  +1\right)^{1/2}\right)\right) {dy \over \sqrt y}. $$
Then

$$\left| \int_{1}^{(1+\tau)^{2\delta}}\sin \left( 2xy -  {\pi\over 4} \right) \cos \left(2 \tau \log \left( y + \left( y^2  +1\right)^{1/2}\right)\right) {dy \over \sqrt y}\right|$$

$$\le \int_{1}^{(1+\tau)^{2\delta}} {dy \over \sqrt y} = 2 \left[ (1+\tau)^\delta -1\right].$$
Concerning the second integral, we have,  owing to the integration by parts, 

$$ \int_{(1+\tau)^{2\delta}}^\infty  \sin \left( 2xy -  {\pi\over 4} \right) \cos \left(2 \tau \log \left( y + \left( y^2  +1\right)^{1/2}\right)\right) {dy \over \sqrt y} $$

$$= {(1+\tau)^{- \delta}\over 2x}  \cos \left( 2x (1+\tau)^{2\delta} -  {\pi\over 4} \right) \cos \left(2 \tau \log \left( (1+\tau)^{2\delta} + \left( (1+\tau)^{4\delta}  +1\right)^{1/2}\right)\right)$$

$$- {\tau\over x}  \int_{(1+\tau)^{2\delta}}^\infty  \cos \left( 2xy -  {\pi\over 4} \right) \sin \left(2 \tau \log \left( y + \left( y^2  +1\right)^{1/2}\right)\right) {dy \over \sqrt {y (y^2+1)} } $$

$$- {1\over 4x }  \int_{(1+\tau)^{2\delta}}^\infty  \cos \left( 2xy -  {\pi\over 4} \right) \cos \left(2 \tau \log \left( y + \left( y^2  +1\right)^{1/2}\right)\right) {dy \over y^{3/2}},$$
and, accordingly,

$$\left|  \int_{(1+\tau)^{2\delta}}^\infty  \sin \left( 2xy -  {\pi\over 4} \right) \cos \left(2 \tau \log \left( y + \left( y^2  +1\right)^{1/2}\right)\right) {dy \over \sqrt y} \right| $$

$$\le {1\over x}  \ (1+2\tau) (1+\tau)^{- \delta}.$$
Thus, combing with (1.25), we establish the following inequality for the Kontorovich-Lebedev kernel $\left(x, \tau > 0, \ \delta > 0,\  M, N \in \mathbb{N}\right)$

$$    \left|K_{i\tau}(x)\right| \le {\sqrt \pi\over \sqrt{ \tau \sinh(\pi\tau)}} \bigg[ {2\over \sqrt {\pi x} }  \ \left(1 +2\tau \right) (1+\tau)^{- \delta} $$

$$ + {4 \sqrt x\over \sqrt \pi}  \left[ (1+\tau)^\delta -1\right] +1 +  2 b \sqrt{2x}  \left( 1 + \alpha + {3 \alpha^2\over 10}\right)^{1/2} $$

$$+  {4\sqrt x \over \pi } \left[ \sum_{n=1}^{N-1}  2^{2n-1} \Gamma \left(2n-{1\over 2}\right) B\left(2n+{3\over 2}, 2n- {1\over 2}\right) \right.$$

$$\left.  +  \sum_{n=0}^{M-1}   4^n  \Gamma\left(2n+{1\over 2}\right) \ B\left(2n+{5\over 2}, 2n+{1\over 2}\right)\right]$$

$$+   {2\sqrt {x} \over \pi^2 }\left[ 4^N \  \Gamma \left(2N-{3\over 2}\right) B\left(2N-{3\over 2}, \ 2N+{1\over 2}\right) \right. $$

$$\bigg. \left. +   4^{M-1} \Gamma \left(2M-{1\over 2}\right) B\left(2M-{1\over 2}, \ 2M+{3\over 2}\right)\right] \bigg]^{1/2}.\eqno(1.26) $$

Finally in this section we show that the integral (see [3], Vol. II, Entry 2.16.3.6)

$$K^2_{i\tau}(x)= 2 \int_1^\infty {K_{2i\tau}(2xy) \over (y^2-1)^{1/2}} \ dy\eqno(1.27)$$
provides a set of the so-called iterated Lebedev type upper bounds. In fact, owing to (1.5), we find from (1.27)

$$ K^2_{i\tau}(x)\le    {(2/x)^{1/4}\   \Gamma(1/4) \over \sqrt{\sinh (2\pi\tau)}} \int_1^\infty {y^{-1/4} \over (y^2-1)^{1/2}} \ dy$$

$$=  { 2^{- 3/2}\  \Gamma^2(1/8)  \  x^{-1/4} \over \sqrt{\sinh (2\pi\tau)}},$$
i.e. we derive the following inequality

$$ \left| K_{i\tau}(x) \right| \le      {\Gamma(1/8)\over  2^{3/4}}   \ { x^{-1/8} \over \sinh^{1/4} (2\pi\tau)},\quad x, \tau > 0.\eqno(1.28)$$
Now, returning to (1.27), we apply (1.28) to get the estimate

$$ K^2_{i\tau}(x) \le      \ { (2/x)^{1/8}\ \Gamma(1/8)  \over \sinh^{1/4} (4\pi\tau)} \int_1^\infty {y^{-1/8} \over (y^2-1)^{1/2}} \ dy = \ {   2^{-7/4} \Gamma^2(1/16) \  x^{-1/8} \over \sinh^{1/4} (4\pi\tau)}.$$
Hence it yields the inequality

$$ \left| K_{i\tau}(x) \right| \le   {    \Gamma(1/16)\over 2^{7/8}} \  {x^{-1/16} \over \sinh^{1/8} (4\pi\tau)}, \quad x, \tau > 0.\eqno(1.29)$$
Continuing to apply the bound which is obtained on the $n$-th step to integral (1.27), we observe the sequence of the Lebedev type inequalities 

$$ \left| K_{i\tau}(x) \right| \le   { \Gamma\left(2^{-n-1}\right)\over  2^{1- 2^{-n}} } \  \bigg[  \sqrt  x \sinh \left(2^{n} \pi\tau /2\right)\bigg]^{- 2^{-n}},  \quad n \in \mathbb{N},\quad x, \tau > 0.\eqno(1.30)$$

\section {Uniform asymptotic expansion}

In this section we propose a new approach to establish a uniform asymptotic expansion for the kernel $K_{i\tau}(x)$, comparing with [1, Section 1.2],  and give  an explicit remainder term  with the corresponding error estimate.   We start, appealing to its relation with the modified Bessel function of the first kind

$$K_{i\tau}(x)= {\pi\over 2i\sinh(\pi\tau)} \left[ I_{-i\tau}(x)- I_{i\tau}(x)\right],\quad x,  \tau  >0,\eqno(2.1)$$ 
where

$$I_{\nu}(x)= \sum_{k=0}^\infty {(x/2)^{2k+\nu} \over k! \Gamma(k+\nu+1)}.\eqno(2.2)$$
Hence, employing the reflection formula for gamma function, we write (2.1) in the form

$$ K_{i\tau}(x)= {\Gamma(i\tau) \over 2}  \sum_{k=0}^\infty {(x/2)^{2k-i\tau} \Gamma (1-i\tau) \over k! \Gamma(k-i\tau+1)}  +  {\Gamma(-i\tau) \over 2}  \sum_{k=0}^\infty {(x/2)^{2k+i\tau}  \Gamma (1+i\tau) \over k! \Gamma(k+i\tau+1)}$$

$$=  {\rm Re} \left[ \Gamma(i\tau) \left({x\over 2}\right)^{-i\tau} \right]  + {\Gamma(i\tau) \over 2}  \sum_{k=1}^\infty {(x/2)^{2k-i\tau} \Gamma (1-i\tau) \over k! \Gamma(k-i\tau+1)} $$

$$+  {\Gamma(-i\tau) \over 2}  \sum_{k=1}^\infty {(x/2)^{2k+i\tau}  \Gamma (1+i\tau) \over k! \Gamma(k+i\tau+1)}.$$  
Then, appealing to the simple beta-integral, we obtain the equality

$$ K_{i\tau}(x)=  {\rm Re} \left[ \Gamma(i\tau) \left({x\over 2}\right)^{-i\tau} \right]  + {\Gamma(i\tau) \over 2}  \sum_{k=0}^\infty {(x/2)^{2k-i\tau+2}  \over (k+1)!\  k!} \int_0^1 y^{-i\tau}(1-y)^k dy $$

$$+  {\Gamma(-i\tau) \over 2}  \sum_{k=0}^\infty {(x/2)^{2k+i\tau+2}  \over (k+1)!\ k!} \int_0^1 y^{i\tau}(1-y)^k dy.\eqno(2.3)$$  
Interchanging the order of integration and summation via the absolute and uniform convergence and using (2.2), we derive after simple substitutions

$$ K_{i\tau}(x)=  {\rm Re} \bigg[ \Gamma(i\tau) \left({x\over 2}\right)^{-i\tau} \left[ 1  +  x^{2i\tau}  \int_0^x (x^2-y^2)^{-i\tau} I_1\left(  y \right) dy \right] \bigg].\eqno(2.4)$$  
The integral in (2.4) can be treated via integration by parts and the use of the differential relation for the modified Bessel function of the first kind $I_\nu(z)$ [3, Vol. II]

$${d\over dz} \left[ z^{-\nu} I_\nu (z) \right] = z^{-\nu} I_{\nu+1}(z).\eqno(2.5)$$
Indeed, we find after   $N$ times consecutive  integration by parts

$$   x^{2i\tau} \int_0^x (x^2-y^2)^{-i\tau} I_1\left(  y \right) dy = {x^2\over 4 (1-i\tau)} +  { x^{2i\tau} \over 2(1-i\tau)}  \int_0^x (x^2-y^2)^{1-i\tau} I_2\left(  y \right) {dy\over y}$$

$$= \sum_{m=1}^{N} {(x/2)^{2m}\over   m! \ (1-i\tau)_m} + { x^{2i\tau} \over 2^N   (1-i\tau)_N}  \int_0^x (x^2-y^2)^{N-i\tau} I_{N+1}\left(  y \right) {dy\over y^N},\eqno(2.6)$$
where $(a)_n$ is the Pochhammer symbol.  Thus we derive the following key formula for the asymptotic expansion of the Kontorovich-Lebedev kernel by  the index 

$$ K_{i\tau}(x)=  {\rm Re} \bigg[ \Gamma(i\tau) \left({x\over 2}\right)^{-i\tau} \left[ 1 +  \sum_{m=1}^{N} {(x/2)^{2m}\over 2^m m! \ (1-i\tau)_m}\right.$$

$$\left. + { x^{2i\tau} \over 2^N  \ (1-i\tau)_N}  \int_0^x (x^2-y^2)^{N-i\tau} I_{N+1}\left(  y \right) {dy\over y^N}\right] \bigg],\quad x, \tau > 0,\ N \in \mathbb{N}_0.\eqno(2.7)$$
Precisely, we have

{\bf Theorem 1}. {\it Let $N$ be a non-negative  integer and $ x \in \left(0, X\right],\ X >0$. Then the modified Bessel function $K_{i\tau}(x)$  has the following asymptotic expansion

$$ K_{i\tau}(x)=  \sqrt{{2\pi\over \tau}}\  e^{-\pi\tau/2} \bigg[ \cos \left( \tau \log \left({2\tau\over e x}\right)  - {\pi\over 4} \right) + R_N (\tau) \bigg],\quad \tau \to +\infty,\eqno(2.8)$$
where the remainder term is given explicitly

$$ R_N (\tau) =  {\rm Re} \bigg[ \exp \left( i \left( \tau \log \left({2\tau\over e x}\right) - {\pi\over 4} \right) \right) \bigg[ r(\tau) + \bigg( 1+ r(\tau) \bigg)  \bigg.\bigg.$$

$$\left.\left.  \times \left[ \sum_{m=1}^{N} {(x/2)^{2m}\over  m! \ (1-i\tau)_m} + { x^{2i\tau} \over 2^N  (1-i\tau)_N}  \int_0^x (x^2-y^2)^{N-i\tau} I_{N+1}\left(  y \right) {dy\over y^N}\right] \right] \right],\eqno(2.9)$$
and 

$$r(\tau) = \exp\left( \int_0^\infty e^{-i\tau t} \left[ {1\over 2}- {1\over t} + {1\over e^t-1}\right] {dt\over t} \right) -1.\eqno(2.10)$$
Moreover, the remainder term $R_N(\tau)$ has the following upper bound

$$  \left| R_N (\tau) \right| \le  {1\over \tau} \left[ {e^{1/(6\tau_0)} \over 6}+  \left( \tau_0 +  { e^{1/(6\tau_0)} \over 6 }\right)  \left[ \exp \left( {X^2\over 4\tau_0} \right) \right.\right.$$

$$\left.\left.+ \left( {X^2 \over 2\tau_0}\right)^N  \left[ {I_{N}\left( X \right) \over X^N} - {1\over 2^N N!} \right]\right] \right],\eqno(2.11)$$
where $\tau \ge \tau_0 > 0$.}

\begin{proof} In fact, appealing to the Stirling formula for gamma function [5] for the pure imaginary argument, we have the equality

$$\Gamma(i\tau)   \left({x\over 2}\right)^{-i\tau} =  \sqrt{{2\pi\over \tau}}\  \exp \left(- {\pi\tau\over 2}+ i \left( \tau \log \left({2\tau\over e x}\right)  - {\pi\over 4} \right)\right)
\bigg[ 1+ r(\tau)\bigg],\eqno(2.12)$$
where $r(\tau)$ is defined by (2.10).  Then, substituting the right-hand side of (2.11) in (2.7) and taking the real part, we get (2.8).  In order to prove inequality (2.11), we use the known bound for the remainder $r(\tau)$ [5], namely,

$$\left|r(\tau)\right| \le e^{1/(6\tau)} -1,\eqno(2.13)$$
where $\tau $ is bounded away from zero, i.e. $\tau \ge \tau_0 > 0$.  Hence,  from (2.13) we derive the estimate

$$\left|r(\tau)\right| \le {1\over 6\tau} \sum_{k=0}^\infty { (6\tau_0)^{-k}\over (k+1)!} \le  { e^{1/(6\tau_0)} \over 6\tau}.$$
Therefore, we have from (2.2), (2.5), (2.9)

$$ \left| R_N (\tau) \right| \le \left| r(\tau) \right| + \bigg( 1+ \left| r(\tau) \right| \bigg)$$

$$ \times \left[ \sum_{m=1}^{N} {(X/2)^{2m}\over  m! \ |(1-i\tau)_m|} + { 1 \over 2^N  |(1-i\tau)_N|}  \int_0^X (X^2-y^2)^{N} I_{N+1}\left(  y \right) {dy\over y^N} \right]$$

$$\le {1\over \tau} \left[{e^{1/(6\tau_0)} \over 6}+  \left( \tau_0 +  { e^{1/(6\tau_0)} \over 6} \right)\left[ \sum_{m=1}^{N} {(X/(2\sqrt {\tau_0}) )^{2m}\over  m!} + {1\over (2\tau_0)^N }  \int_0^X (X^2-y^2)^{N} I_{N+1}\left(  y \right) {dy\over y^N} \right]\right]$$

$$\le {1\over \tau} \left[{e^{1/(6\tau_0)} \over 6}+  \left( \tau_0 +  { e^{1/(6\tau_0)} \over 6 }\right)  \left[ \exp \left( {X^2\over 4\tau_0} \right) + \left( {X^2 \over 2\tau_0}\right)^N  \left[ {I_{N}\left( X \right) \over X^N} - {1\over 2^N N!} \right]\right] \right].$$
This proves (2.11) and completes the proof of Theorem 1. 

\end{proof}

\section{Summability of the KL-integrals in the Jones sense}

As it is suggested in [6, Section 5], we will consider the following Kontorovich-Lebedev integrals in the Jones sense [7]

$$f(x,a) \equiv  \lim_{\varepsilon \to 0+} f_\varepsilon(x,a)=  \lim_{\varepsilon \to 0+} \int_0^\infty e^{-\varepsilon \tau^2} \left[ \psi_1(\tau) \cosh\left({\pi\over 2}  +a \right)\tau\right.$$

$$\left.+ \psi_2(\tau)\tau \sinh\left({\pi\over 2} +a \right)\tau \right] K_{i\tau} (x) d\tau,\eqno(3.1)$$
where $x >0$, $a$ is a parameter  and $\psi_1(\tau),\ \psi_2(\tau) $ are even entire functions of the exponential type.  Our goal is to prove the existence of the limit (3.1) in a weak sense under the restriction ${\rm Re}  a \in [0, \pi/2)$.  First we observe that when $-\pi < {\rm Re} a < 0$ the limit (3.1) exists, appealing to  Theorem 1 and Lebedev's type upper bounds above via the absolute and uniform convergence.  Second, it is enough to consider the real  case of $a \in [0, \pi/2)$.

We begin with a key example of (3.1).

{\bf Theorem 2}.  {\it Let $x > 0,\ a \in [0,\pi/2)$. The  following limit holds 

$$ \lim_{\varepsilon \to 0+} \int_0^\infty e^{-\varepsilon \tau^2}  \cosh \left[\left({\pi\over 2}  +a \right)\tau\right] K_{i\tau} (x) d\tau $$

$$\equiv   \lim_{\varepsilon \to 0+} I_\varepsilon (x,a)= {\pi\over 2} \exp\left(x\sin a\right),\eqno(3.2)$$
where the convergence is understood in the sense of the generalized Mellin transform (cf. $[8]$), i.e.}

$$\lim_{\varepsilon \to 0+} \langle  I_\varepsilon (x,a),\  e^{-x} x^{s-1} \rangle = {\pi\over 2} \langle  e^{x\sin a},\  e^{-x} x^{s-1} \rangle,\quad {\rm Re} s > 0.\eqno(3.3)$$

\begin{proof} In fact,  since $e^{-x} I_\varepsilon(x,a) x^{s-1} \in L_1(\mathbb{R}_+)$,  $e^{-x} I_\varepsilon$ represents a regular generalized function. Hence we write, using the conventional Mellin transform, 

$$\int_0^\infty e^{-x} I_\varepsilon(x,a) x^{s-1} dx = \int_0^\infty e^{-x}  x^{s-1} \int_0^\infty e^{-\varepsilon \tau^2}  \cosh \left[\left({\pi\over 2}  +a \right)\tau\right]$$

$$\times  K_{i\tau} (x) d\tau dx.\eqno(3.4)$$ 
The interchange of the order of integration in (3.4) is allowed for each $\varepsilon > 0$ due to the absolute convergence of the iterated integral.  Then, appealing to the   Mellin transform  (cf. Entry 8.4.23.3 in [3, Vol. III])

$$ \int_0^\infty e^{-x} K_{i\tau}(x) x^{s-1} dx =  2^{-s} \sqrt \pi\  \frac{\Gamma(s+i\tau) \Gamma(s-i\tau)}{\Gamma(s+1/2)},\quad {\rm Re} s > 0,\eqno(3.5)$$
we obtain 
$$\int_0^\infty e^{-x} I_\varepsilon(x) x^{s-1} dx =    { 2^{-s} \sqrt \pi\over  \Gamma(s+1/2)}   \int_0^\infty  e^{-\varepsilon \tau^2}  \cosh \left[\left({\pi\over 2}  +a \right)\tau\right]$$

$$\times  \Gamma(s+i\tau) \Gamma(s-i\tau) d\tau.\eqno(3.6)$$
Now, one can pass to the limit when $\varepsilon \to 0+$ under the integral sign on the right-hand side of (3.6) via Abel's test and Stirling's asymptotic formula for  the 
gamma function.  Hence the integral by $\tau$ is calculated in [1, formula (1.104)], namely, 

$$ \int_0^\infty  \cosh \left[\left({\pi\over 2}  +a \right)\tau\right] \Gamma(s+i\tau) \Gamma(s-i\tau) d\tau = {\pi \Gamma(2s) \over 2^{2s}} \left[\cos \left( {\pi\over 4}  + {a\over 2} \right)\right]^{-2s}.\eqno(3.7)$$
Therefore,  employing the duplication formula for gamma function, we find from (3.6), (3.7)

$$\lim_{\varepsilon \to 0+}   \int_0^\infty e^{-x} I_\varepsilon(x,a) x^{s-1} dx = {\pi\over 2} \ \Gamma(s) \left[1- \sin a\right]^{-s}. $$
Finally, employing the Euler integral for the gamma function, we establish (3.2) and complete the proof of Theorem 2.  

\end{proof}

Now we are ready to prove the main result of this section.

{\bf Theorem 3}. {\it  Let $x > 0,\ a \in [0,\pi/2)$ and $\psi_j(\tau), j=1,2$ are even entire functions of the exponential type $b_j \in [0,  (1-\sin a)/(2e) ),\ j=1,2$.  Then the limit $(3.1)$ has the value

$$f (x,a)= {\pi\over 2} \bigg[ \psi_1( D_a) + D_a \psi_2(D_a) \bigg] \left\{ e^{x\sin a} \right\},\eqno(3.8)$$
where $\psi_j(D_a),\ j=1,2$ are functions of the differential operator $D_a\equiv d/da$, acting on $e^{x\sin a}$ and the convergence is understood in the sense of the generalized Mellin transform }

$$\lim_{\varepsilon \to 0+} \langle  f_\varepsilon (x,a),\  e^{-x} x^{s-1} \rangle $$

$$=  {\pi\over 2} \langle  \left[ \psi_1( D_a) + D_a \psi_2(D_a) \right] \left\{ e^{x\sin a} \right\} ,\  e^{-x} x^{s-1} \rangle,\  {\rm Re} s > 0.\eqno(3.9)$$

\begin{proof} Indeed, functions $\psi_j(\tau),\ j=1,2$ are represented in terms of the Taylor series

$$\psi_j(\tau) = \sum_{n=0}^\infty c_{2n,j} \tau^{2n},\quad  j=1,2,$$
where the coefficients $c_{2n,j} $ satisfy the Cauchy estimates 

$$|c_{2n,j} | <  \left( { e b_j\over 2n}\right)^{2n},\quad  n > N,\ j=1,2.\eqno(3.10)$$
Hence by the same arguments as in the proof of Theorem 2, we recall (3.1), (3.5) to obtain

$$\int_0^\infty e^{-x} f_\varepsilon(x,a) x^{s-1} dx =  { 2^{-s} \sqrt \pi\over  \Gamma(s+1/2)}  \int_0^\infty e^{-\varepsilon \tau^2} \left[ \psi_1(\tau) \cosh \left[\left({\pi\over 2}  +a \right)\tau\right] \right.$$

$$\left.+ \psi_2(\tau)\tau \sinh \left[\left({\pi\over 2} +a \right)\tau \right] \right]  \Gamma(s+i\tau) \Gamma(s-i\tau)  d\tau.\eqno(3.11)$$
The interchange of the order of integration on the left-hand side of is guaranteed  for each positive $\varepsilon $ via the absolute convergence of the corresponding iterated integral.  Then in the same manner one passes to the limit when $\varepsilon \to 0+$ due to the convergence of the integral

$$  \int_0^\infty  \left[ \psi_1(\tau) \cosh \left[ \left({\pi\over 2}  +a \right)\tau \right] + \psi_2(\tau)\tau \sinh\left[ \left({\pi\over 2} +a \right)\tau \right] \right]  \Gamma(s+i\tau) \Gamma(s-i\tau)  d\tau $$
under conditions of the theorem.  Therefore,

$$\lim_{\varepsilon \to 0+} \int_0^\infty e^{-x} f_\varepsilon(x,a) x^{s-1} dx =  { 2^{-s} \sqrt \pi\over  \Gamma(s+1/2)}  \int_0^\infty  \left[ \psi_1(\tau) \cosh \left[\left({\pi\over 2}  +a \right)\tau \right]\right.$$

$$\left.+ \psi_2(\tau)\tau \sinh \left[ \left({\pi\over 2} +a \right)\tau \right]\right]  \Gamma(s+i\tau) \Gamma(s-i\tau)  d\tau.\eqno(3.12)$$
Meanwhile,  returning to (3.7) and differentiating with respect to parameter $a$, we get the value 

$$ \int_0^\infty  \tau \sinh \left[\left({\pi\over 2}  +a \right)\tau\right] \Gamma(s+i\tau) \Gamma(s-i\tau) d\tau$$

$$ = {\pi \Gamma(2s+1) \over 2^{2s+1}} \left[\cos \left( {\pi\over 4}  + {a\over 2} \right)\right]^{-2s-1} \sin \left( {\pi\over 4}  + {a\over 2} \right),\eqno(3.13)$$
since the latter integral converges uniformly for $0\le a \le \pi/2- a_0,\ a_0 \in (0, \pi/2].$  The right-hand side of (3.12) can be treated by substitution of the series (3.9) and using the differentiation with respect to a parameter.  In fact,  it yields 

$$\sum_{n=0}^\infty c_{2n,1}  \int_0^\infty  \tau^{2n} \cosh \left[ \left({\pi\over 2}  +a \right)\tau \right]  \Gamma(s+i\tau) \Gamma(s-i\tau)  d\tau$$

$$+  \sum_{n=0}^\infty c_{2n,2} \int_0^\infty  \tau^{2n+1}  \sinh \left[\left({\pi\over 2} +a \right)\tau \right]  \Gamma(s+i\tau) \Gamma(s-i\tau)  d\tau,$$
where the interchange of the order of integration and summation and further differentiation under the integral sign can be justified, employing inequality (3.10),  the representation of the product of gamma functions by  the Mellin transform  from Entry  8.4.23.1 in [3, Vol. II]

$$\Gamma(s+i\tau) \Gamma(s-i\tau) = 2^{2(1-s)} \int_0^\infty  K_{2i\tau}(x) x^{2s-1} dx,\quad {\rm Re} s > 0\eqno(3.14)$$
and the inequality for the Macdonald function (see [1, formula (1.100))

$$|K_{i\tau} (x)| \le e^{-\delta \tau}K_0(x\cos\delta ),\quad x, \tau >0,\ \delta \in \left[0, {\pi\over 2} \right).\eqno(3.15)$$
Precisely,  we derive for a big enough positive integer $N$ and $ \delta \in ( 1/2 (a+\pi/2+ \max (b_1,b_2)), \ \pi/2) $

$$\sum_{n=N+1}^\infty \left| c_{2n,1} \right|  \int_0^\infty  \tau^{2n} \cosh \left[ \left({\pi\over 2}  +a \right)\tau \right] \left| \Gamma(s+i\tau) \Gamma(s-i\tau) \right|  d\tau$$

$$+  \sum_{n=N+1}^\infty \left|c_{2n,2} \right| \int_0^\infty  \tau^{2n+1}  \sinh \left[\left({\pi\over 2} +a \right)\tau \right]  \left|\Gamma(s+i\tau) \Gamma(s-i\tau) \right| d\tau$$

$$<  2^{3- 2 {\rm Re} s} \int_0^\infty  K_{0}(x \cos\delta ) x^{2{\rm Re} s -1} dx \left[  \sum_{n=N+1}^\infty \left( { e b_1\over 2n}\right)^{2n}   \int_0^\infty  \tau^{2n} e^{- (2\delta -a -\pi/2 ) \tau} d\tau  \right.$$

$$\left. +  \sum_{n=N+1}^\infty \left( { e b_2\over 2n}\right)^{2n}   \int_0^\infty  \tau^{2n+1} e^{- (2\delta -a -\pi/2 ) \tau} d\tau \right]$$

$$=  {2^{3- 2 {\rm Re} s} \over 2\delta-a-\pi/2} \int_0^\infty  K_{0}(x \cos\delta ) x^{2{\rm Re} s -1} dx  \sum_{n=N+1}^\infty (2n)! \left( { e \over 2n}\right)^{2n} \left[ \left( {  b_1\over 2\delta-a-\pi/2 }\right)^{2n}  \right.$$

$$\left. +   {2n+1\over 2\delta-a-\pi/2}  \left( {  b_2\over 2\delta-a-\pi/2}\right)^{2n}   \right] \to 0,\quad N \to \infty.$$
Hence we find from (3.7), (3.12), (3.13) the equalities

$$ \lim_{\varepsilon \to 0+} \int_0^\infty e^{-x} f_\varepsilon(x,a) x^{s-1} dx =  { 2^{-s} \sqrt \pi\over  \Gamma(s+1/2)} \left[\sum_{n=0}^\infty c_{2n,1}  {d^{2n}\over d a^{2n}} \int_0^\infty   \cosh \left[ \left({\pi\over 2}  +a \right)\tau \right] \right.$$ 

$$\left. \times \Gamma(s+i\tau) \Gamma(s-i\tau)  d\tau  +  \sum_{n=0}^\infty c_{2n,2}  {d^{2n+1}\over d a^{2n+1}} \int_0^\infty  \cosh \left[\left({\pi\over 2} +a \right)\tau \right]  \Gamma(s+i\tau) \Gamma(s-i\tau)  d\tau\right]$$
$$=  {\pi\over 2} \Gamma(s) \left[\sum_{n=0}^\infty c_{2n,1}  {d^{2n}\over d a^{2n}}   \left[ 1- \sin a \right]^{-s} +  \sum_{n=0}^\infty c_{2n,2}  {d^{2n+1}\over d a^{2n+1}}   \left[1-\sin a \right]^{-s}  \right]$$

$$= {\pi\over 2} \bigg[ \psi_1( D_a) + D_a \psi_2(D_a) \bigg]  \int_0^\infty \exp \left(- x( 1- \sin a )\right) x^{s-1} dx,\quad D_a\equiv {d\over da}.\eqno(3.15)$$ 
The differentiation with respect to $a$ under the integral sign of the latter integral in (3.15) and the action of the operator $\psi_1( D_a) + D_a \psi_2(D_a)$ inside this integral  we will motivate, as above,  owing to the uniform convergence of the integrals for   derivatives. To do this,  we appeal to the Hoppe formula [9] to write the $2n$th derivative of  $\exp (x\sin a) $ in the form

$${d^{2n}\over da^{2n}} \left[ e^{x\sin a}\right] = e^{x\sin a} \sum_{k=0}^{2n} {(-1)^k x^k\over k! }\sum_{j=0}^k (-1)^j \binom{k}{j} \left[ \sin a \right]^{k-j} 
{d^{2n}\over da^{2n}} \left[ \sin a\right]^{j} $$

$$= e^{x\sin a} (-1)^n \sum_{k=0}^{2n} {(-1)^k x^k\over k! }\sum_{j=0}^k {(-1)^j \over (2i)^j} \binom{k}{j} \left[ \sin a \right]^{k-j}  \sum_{r=0}^j (-1)^r \binom{j}{r} e^{i(j- 2r) a} (j-2r)^{2n}. $$
Hence

$$\left| {d^{2n}\over da^{2n}} \left[ e^{x\sin a}\right]  \right| \le e^{x\sin a} \sum_{k=0}^{2n} {x^k\over k! }\sum_{j=0}^k {j^{2n} \over 2^j} \binom{k}{j} \sum_{r=0}^j  \binom{j}{r} \le e^{ x \sin a} (2n)^{2n} \sum_{k=0}^{2n} { (2x)^k\over k! }.$$
Then for ${\rm Re} s = \gamma > 0$ and some $n_0 \in \mathbb{N}$

$$ \int_0^\infty  \sum_{n=n_0+1}^\infty  \left| c_{2n,1} \right| \left| {d^{2n}\over da^{2n}} \left[ e^{x\sin a}\right]  \right|   e^{-x} x^{\gamma -1} dx $$

$$\le (1-\sin a)^{-\gamma} \Gamma(\gamma) \sum_{n=n_0+1}^\infty  \left| c_{2n,1} \right|  (2n)^{2n} \sum_{k=0}^{2n} { 2^k (\gamma)_k \over k! (1- \sin a)^k } $$

$$<  (1- \sin a)^{-\gamma} \Gamma(\gamma) \sum_{n=n_0+1}^\infty  (e b_1)^{2n}  \left[ \sum_{k=0}^{n} { 2^{2k} (\gamma)_{2k} \over (2k)! (1- \sin a)^{2k} } \right.$$

$$\left. + \sum_{k=0}^{n-1} { 2^{2k+1} (\gamma)_{2k+1} \over (2k+1)! (1- \sin a)^{2k+1} } \right] 
<  (1- \sin a_0)^{-\gamma} \Gamma(\gamma) \left[ \sum_{k=0}^{\infty} { 2^{2k} (\gamma)_{2k} \over (2k)! (1- \sin a_0)^{2k} }  \sum_{n= k}^\infty  (e b_1)^{2n}  \right.$$

$$\left. +  \sum_{k=0}^{\infty} { 2^{2k+1} (\gamma)_{2k+1} \over (2k+1)! (1- \sin a_0)^{2k+1} } \sum_{n= k}^\infty  (e b_1)^{2(n+1)}\right]$$

$$= {(1- \sin a_0)^{-\gamma} \Gamma(\gamma)\over  1- (e b_1)^2 } \left[ \sum_{k=0}^{\infty} { (2 e b_1)^{2k} (\gamma)_{2k} \over (2k)! (1- \sin a_0)^{2k} } \right.$$

$$\left.   +  2 (e b_1)^2 \sum_{k=0}^{\infty} { (2 e b_1)^{2k} (\gamma)_{2k+1} \over (2k+1)! (1- \sin a_0)^{2k+1} } \right] < \infty$$
when $ 0 \le b_1 < (1-\sin a_0) /(2e),\ 0 \le a\le a_0 < \pi/2 $.   In the same manner we justify the action of the operator $D_a \psi_2(D_a)$ inside the integral in (3.15). 
Thus we arrive at (3.9) and complete the proof of Theorem 3. 

\end{proof}

\bigskip
\noindent{{\bf Acknowledgments}}
\bigskip

\noindent The work was partially supported by CMUP, which is financed by national funds through FCT (Portugal)  under the project with reference UIDB/00144/2020.

\bigskip
\noindent{{\bf Disclosure statement}}
\bigskip

\noindent No potential conflict  of interest was reported by the author.

\bigskip
\centerline{{\bf References}}
\bigskip
\baselineskip=12pt
\medskip
\begin{enumerate}

\item[{\bf 1.}\ ]  Yakubovich S.  Index Transforms. Singapore:  World Scientific Publishing Company;  1996.

\item[{\bf 2.}\ ]   Olenko AY.  Upper bound on $\sqrt x J_\nu(x)$  and its applications. Integral Transforms and Special Functions.  2006; 17, N 6:  455-467.

\item[{\bf 3.}\ ]  Prudnikov AP,  Brychkov  YuA,  Marichev OI. Integrals and series:  Vol. I: Elementary functions. New York:  Gordon and Breach;   1986;   Vol. II:  Special functions. New York: Gordon and Breach;  1986;   Vol. III:  More special functions. New York:   Gordon and Breach; 1990.

\item[{\bf 4.}\ ]   Nemes  G.  Error bounds for the large-argument asymptotic expansions of the Hankel and Bessel functions.  Acta Appl. Math.  2017; 150:  141-177.

\item[{\bf 5.}\ ]  Lebedev NN.  Special Functions and Their Applications. Englewood Cliffs, N.J. : Prentice-Hall, INC; 1965.

\item[{\bf 6.}\ ]  Ehrenmark U. Summability experiments with a class of divergent inverse Kontorovich-Lebedev transforms.  Comput. Math. Appl. 2018; 76, N 1: 141-154.

\item[{\bf 7.}\ ]  Jones DS.  The Kontorovich-Lebedev transform. J. Inst. Math. Appl. 1980; 26: 133-141.  

\item[{\bf 8.}\ ]  Zemanian, A. H. Generalized integral transformations. New York:  Dover Publications; 1987.

\item[{\bf 9.}\ ]  Johnson, WP.  The curious history of Fa{\' a} di Bruno's formula.  Amer. Math. Monthly. 2002; 109, N 3: 217-234.

\end{enumerate}

\end{document}